
\input gtmacros
\input amsnames
\input amstex
\input rlepsf        
\input newinsert   
       
\catcode`\@=12

\input gtmonout
\volumenumber{2}
\volumeyear{1999}
\volumename{Proceedings of the Kirbyfest}
\pagenumbers{299}{320}
\papernumber{16}
\received{7 December 1998}\revised{23 May 1999}
\published{20 November 1999}

\let\\\par
\def\topmatter{\relax}

\let\gttitle\title
\def\title#1\endtitle{\gttitle{#1}}
\let\gtauthor\author
\def\author#1\endauthor{\gtauthor{#1}}
\let\gtaddress\address
\def\address#1\endaddress{\gtaddress{#1}}
\let\gtemail\email
\def\email#1\endemail{\gtemail{#1}}
\def\subjclass#1\endsubjclass{\primaryclass{#1}}
\let\gtkeywords\keywords
\def\keywords#1\endkeywords{\gtkeywords{#1}}
\def\heading#1\endheading{{\def\S##1{\relax}\def\\{\relax\ignorespaces}
    \section{#1}}}
\def\head#1\endhead{\heading#1\endheading}

\def\subhead#1\endsubhead{\sh{#1}}
\def\subsubhead#1\endsubsubhead{\sh{#1}}
\def\specialhead#1\endspecialhead{\sh{#1}}
\def\demo#1{\rk{#1}\ignorespaces}
\def\enddemo{\qed\ppar}

\def\qed{\ifmmode\quad\sq\else\hbox{}\hfill$\sq$\par\goodbreak\rm\fi}  
\def\proclaim#1{\rk{#1}\sl\ignorespaces}
\def\endproclaim{\rm\ppar}
\def\cite#1{[#1]}
\newcount\itemnumber
\def\roster{\items\itemnumber=1}
\def\endroster{\enditems}
\let\itemold\item
\def\item{\itemold{{\rm(\number\itemnumber)}}%
\global\advance\itemnumber by 1\ignorespaces}
\def\S{section~\ignorespaces}  
\def\date#1\enddate{\relax}
\def\thanks#1\endthanks{\relax}   
\def\dedicatory#1\enddedicatory{\relax}  
\let\footnote\plainfootnote

\def\Refs{\ppar{\large\bf References}\ppar\bgroup\leftskip=25pt
\frenchspacing\parskip=3pt plus2pt\small}       
\def\endRefs{\egroup}
\def\widestnumber#1#2{\relax}
\def\endrefitem{}
\def\refdef#1#2#3{\def#1{\leavevmode\unskip\endrefitem#2\def\endrefitem{#3}}}
\def\ref{\par}
\def\endref{\endrefitem\par\def\endrefitem{}}
\refdef\key{\noindent\llap\bgroup[}{]\ \ \egroup}
\refdef\no{\noindent\llap\bgroup[}{]\ \ \egroup}
\refdef\by{\bf}{\rm, }
\refdef\manyby{\bf}{\rm, }
\refdef\paper{\it}{\rm, }
\refdef\book{\it}{\rm, }
\refdef\jour{}{ }
\refdef\vol{}{ }
\refdef\yr{$(}{)$ }
\refdef\ed{(}{ Editor) }
\refdef\publ{}{ }
\refdef\inbook{from: ``}{'', }
\refdef\pages{}{ }
\refdef\page{}{ }
\refdef\paperinfo{}{ }
\refdef\bookinfo{}{ }
\refdef\publaddr{}{ }
\refdef\eds{(}{ Editors)}
\refdef\bysame{\hbox to 3 em{\hrulefill}\thinspace,}{ }
\refdef\toappear{(to appear)}{ }
\refdef\issue{no.\ }{ }

\topmatter
\title Simplicial moves on complexes and manifolds
\endtitle

\author W\kern.16em B\kern.16em R Lickorish \endauthor
\asciiauthors{W B R Lickorish}
\address  Department of Pure Mathematics and Mathematical Statistics,
University of Cambridge\\16, Mill Lane, Cambridge, CB2 1SB, UK \endaddress

\email wbrl@dpmms.cam.ac.uk\endemail

\keywords  Simplicial complexes, subdivisions, stellar subdivisions, stellar
manifolds, Pachner moves
\endkeywords

\primaryclass{57Q15}\secondaryclass{52B70}

\abstract 
Here are versions of the proofs of two classic theorems of
combinatorial topology.  The first is the result that piecewise
linearly homeomorphic simplicial complexes are related by stellar
moves.  This is used in the proof, modelled on that of Pachner, of the
second theorem. This states that moves from only a finite collection
are needed to relate two triangulations of a piecewise linear
manifold.
\endabstract
\makeshorttitle

\cl{\small\it For Rob Kirby, a sixtieth birthday offering  
after thirty years of friendship}

\head Introduction \endhead

A finite simplicial complex can be viewed as a combinatorial abstraction or as
a structured subspace of Euclidean space.  The combinatorial approach can lead
to cleaner statements and proofs of theorems, yet it is the more topological
interpretation that provides motivation and application.  In particular,
continuous functions can be approximated by piecewise linear maps, those
generalisations of the simple idea represented by a saw-tooth shaped graph. 
To discuss piecewise linear functions it is necessary to consider
simplex-preserving functions between arbitrary subdivisions of complexes
(subdivisions obtained by disecting Euclidean simplexes into smaller simplexes
in a linear, but unsystematic, way).  The commonly used standard treatments of
piecewise linear topology and piecewise linear manifold theory, those by 
J\,F\,P Hudson \cite 6, C\,P Rourke and B\,J Sanderson \cite {13}  
and E\,C Zeeman \cite {16}, are
successfully based on this idea of arbitrary subdivision.  Indeed, the
newly-found ability to proceed without
 heavy combinatorics was a key factor encouraging the renaissance of piecewise linear topology in the 1960's.
However, one of the main
triumphs of combinatorial topology, as developed by J\,W Alexander  \cite 1  and M\,H\,A Newman \cite 9  in the 1920's
and 1930's, was a result connecting the abstract and piecewise linear approaches.  A proof of that result is given here.   It states
that any subdivision of a finite simplicial complex can be obtained from the original complex by a finite sequence of stellar
moves.  These moves, explained below, change sub-collections of the simplexes in a specified manner, but there are, nonetheless,
infinitely many possible such moves.  A version of this result, not discussed in \cite 6, 
\cite {13} and \cite {16}, does appear in
 \cite 5, but this is not easily available.  Otherwise the original accounts, which use some
outmoded terminology, must be consulted.  The version of stellar theory given here is based on some
lectures given by Zeeman in 1961.  

Inaccessibility of a result hardly matters if nobody wishes to use and understand it.
However in much more recent years, U Pachner  \cite {11}, \cite {12} has, {\it starting} with the
theory of stellar moves, produced a {\it finite} set of combinatorial moves (analogous, in some
way, to the well known Reidemeister moves of knot theory) that suffice to change from any
triangulation of a piecewise linear manifold to any other.  Much of this was foreshadowed by Newman
in \cite 9. A discussion and proof of this, based on Pachner's work, is the task of the final
section of this paper.   As was pointed out by V\,G Turaev and O Ya Viro, such a finite set of moves
can at once be used to establish the well-definedness of any manifold invariant defined in terms of
simplexes.  This they did, for 3--manifolds, in establishing their famous state sum invariants 
\cite {15}  that turned out to be the squares of the moduli of the Witten invariants.  It is not so
very hard to prove the  combinatorial result in dimension three, though still assuming stellar
theory, (see \cite {15} or \cite 2).  However, the result is now being used to substantiate
putative topological quantum field theories in four or more dimensions.  It thus seems desirable to
have available,  readily accessible to topologists, a complete reworking of the stellar theory and
Pachner's work in a common notation as now used by topologists. Aside from their uses, these two
results form a most elegant chapter in the theory of combinatorial manifolds.  In addition, as
explained in \cite 8, the moves between triangulations of a manifold are a simplicial version of
the moves on handle structures of a manifold that come from Cerf theory.   Cerf theory, of course,
has many important uses, one being in Kirby's proof of the sufficiency of the surgery moves of his 
calculus for 3--manifolds.  A piecewise linear version of Cerf theory was not available when any of
the classic piecewise linear texts was written.


\head Notation and other preliminaries \endhead

Some standard notational conventions that will be used will now be outlined, but it is hoped that
basic ideas of simplicial complexes will be familiar. Simplicial
complexes are objects that consist of  finitely many simplexes, which can be thought of as
elementary building blocks, glued together to make up the complex.  Such a complex can be thought of
as just a finite set (the vertices) with certain specified subsets (the simplexes).  Pachner's
results and the basic stellar subdivision theory are probably best viewed in terms of this
latter abstract formulation. Although the abstraction  has its virtues, to be of use in
topology  a simplicial complex needs to represent a topological space. Here that is taken to be a
subspace of some
${\Bbb R}^N$. Thus in what follows a simplex will be taken to be the convex hull of its vertices, they being
independent points in ${\Bbb R}^N$. This allows the use of certain topological words (like `closure' and
`interior') but, more importantly, it allows the idea of an arbitrary subdivision of a complex which, in
turn, permits at once the definition of the natural and ubiquitous idea of a piecewise linear map.  
One of the purposes of this account is to discuss the relation between the  piecewise linear
and  abstract notions, so it may be  wise to feel familiar with both interpretations of a simplicial
complex.  
  The notation
$A
\leq B$ will mean that a simplex
$A$ is a face of simplex $B$; the empty simplex $\emptyset$ is a face of every simplex.

\proclaim {Definition 2.1}\rm A (finite) simplicial  complex $K$ is a finite collection of simplexes, contained
(linearly) in some ${\Bbb R}^N$, such that 
\roster
\item $B \in K$ and $A \leq B$ implies that $A \in K$,
\item  $A \in K$ and $B \in K$ implies that $A \cap B$ is a face of both $A$ and
$B$.
\endroster
\endproclaim

The standard complex $\Delta^n$  will be the complex consisting of all faces of an
$n$--simplex (including the $n$--simplex itself);  its boundary, denoted $\partial \Delta^n$, will be
the subcomplex of all the  {\it proper} faces of $\Delta^n$.  
Sometimes the symbol $A$ will be used ambiguously to denote a simplex $A$ and also
the simplicial complex consisting of $A$ and all its faces.  

The join of simplexes $A$ and $B$ will be
denoted $A\star B$, this being meaningful only when all the vertices of  $A$ and $B$ are
independent. Observe that $\emptyset\star A = A$. The join of simplicial complexes $K$ and $L$,  written $K\star L$,   is
$\{ A\star B :  A \in K , B \in L \} $, where it is assumed that, for $A \in K$ and $B \in L$, the vertices of  $A$ and $B$ are
independent.  Note that the join notation is associative
and commutative.

The link of  a simplex  $A$ in a simplicial complex
$K$, denoted  $\hbox{\rm lk}\, (A,K)$,  is defined by $$\hbox{{\rm lk}}\, (A,K) = \{ B \in K : A\star B \in K
\}.$$   The (closed) star of $A$ in $K$, $\hbox{\rm st}\, (A,K)$, is the join $A\star \hbox{\rm lk}\,
(A,K)$.

A { simplicial isomorphism} between two complexes is a 
bijection between their vertices that induces a bijection between their simplexes.
The {\it polyhedron}
$|K|$ underlying the simplicial complex $K$ is defined to be $\cup_{A \in K}A$ and
$K$ is called a {triangulation} of $|K|$.

A simplicial complex $K^\prime$ is a subdivision of the simplicial
complex $K$ if $|K^\prime | = |K|$ and each simplex of $K^\prime$ is contained linearly in some
simplex of $K$.  Two simplicial complexes $K$ and $L$ are piecewise linearly homeomorphic if they
have subdivisions $K^\prime$ and $L^\prime$ that are simplicially isomorphic.  It is straightforward to show that
this is an equivalence relation.  (More generally, a
piecewise linear map from $K$ to $L$ is  a simplicial map from some subdivision of $K$ to some
subdivision of $L$.) 

\proclaim {Definition 2.2}\rm A combinatorial $n$--ball is a simplicial complex $B^n$ piecewise linearly
homeomorphic to
$\Delta ^n$. A combinatorial $n$--sphere is a simplicial complex $S^n$ piecewise linearly
homeomorphic to $\partial\Delta ^{n+1}$.  A combinatorial $n$--manifold is a simplicial complex $M$
such that, for every vertex $v$ of $M$,
$\hbox {\rm lk}\, (v,M)$ is a combinatorial $(n-1)$--ball or a combinatorial $(n-1)$--sphere. \endproclaim

It could be argued that this traditional definition is not exactly `combinatorial'. The results
described later do show it  to be equivalent to other formulations with a stronger claim to this
epithet.  The definition of a combinatorial $n$--manifold  is easily seen to be equivalent to one
couched in terms of coordinate charts modelled on $\Delta ^n$ with overlap maps being required to
be piecewise linear \cite {13}.  By definition, a piecewise
linear $n$--manifold is just  a class of combinatorial $n$--manifolds equivalent under piecewise linear
homeomorphism.  Beware the fact that a simplicial complex $K$ for which $|K|$ is topologically an
$n$--manifold is not necessarily a combinatorial $n$--manifold. This follows, for example, from the famous theorem of
R\,D Edwards  \cite 4  that states that, if $M$ is a connected orientable combinatorial
3--manifold with the same homology as the 3--sphere (there are infinitely many of these), then 
 $|S^1 \star M|$ is (topologically) homeomorphic to the 5--sphere. However, if $|M|$
is not simply connected, the complex $S^1 \star M$ cannot be a combinatorial 5--manifold.

\proclaim {Definition 2.3}\rm Suppose that $A$  is an $r$--simplex in an abstract simplicial complex 
$K$ and that  $\hbox {\rm lk}\, (A, K) = \partial B$ for some $(n-r)$--simplex
$B \notin K$.      The 
bistellar move $\kappa (A,B)$ consists of changing $K$ by removing
$A\star \partial B$ and inserting $\partial A \star B$.  \endproclaim

\noindent When $n = 2$, the
three types of bistellar move are shown, for $\hbox {\rm dim}\, A = 2, 1, 0$, in \hbox{Figure 1.}
Note that the definition is given for abstract complexes. The condition that $B \notin K$ requires $B$ to be a {\it new} simplex not
seen in $K$. For complexes in ${\Bbb R}^N$ this condition should be replaced by
a requirement that $A \star B$ should exist and that $|A \star B| \cap |K| = |A \star \partial B|$.

\figure\epsfxsize0.5\hsize
\epsfbox{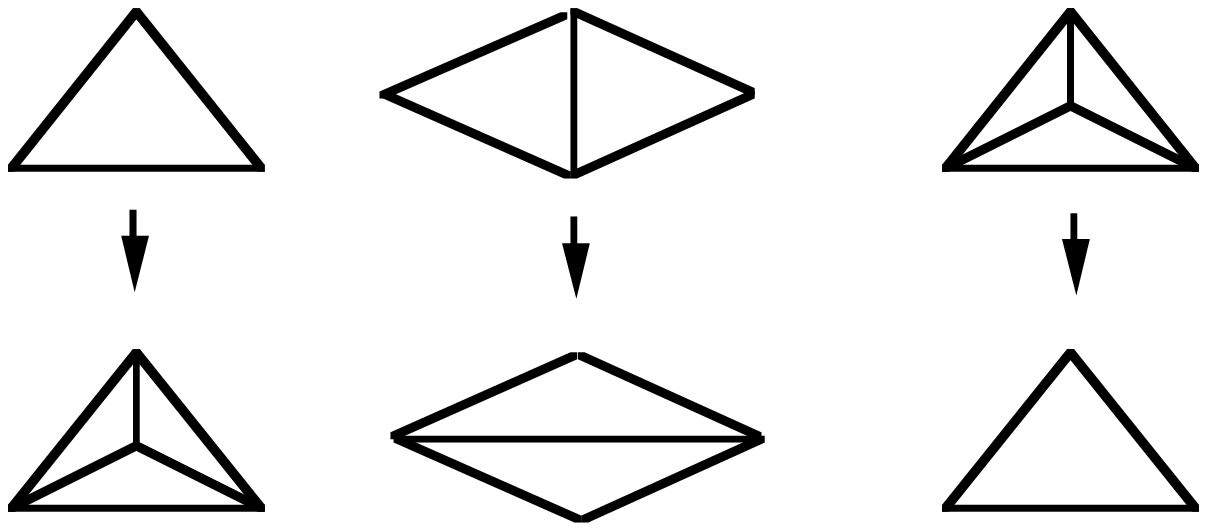}
\endfigure

Complexes related by a finite sequence of these bistellar
moves and simplicial isomorphisms are called {\it bistellar equivalent}.
The main
theorem, described by Pachner \cite {12}, can now be stated as follows.

\proclaim {Theorem 5.9} {\rm(\cite 9, \cite {12})}\qua 
Closed combinatorial $n$--manifolds are piecewise linearly homeomorphic if and
only if they are bistellar equivalent. \endproclaim

\noindent There is also a version of this result for manifolds with boundary that will be discussed at the end of this paper.


\head  Stellar subdivision theory \endhead

Theorem 5.9 is entirely combinatorial in nature.  It can be seen as an extension of
the long known (\cite 1,  \cite 9) combinatorial theory of stellar subdivision.  This stellar
theory will now be described.  In outline, a definition of stellar subdivision leads at once to a
definition of a stellar $n$--manifold, and Pachner's methods can be applied {\it at once} to such
manifolds to prove Theorem 5.9.  The classical theory, described in this section, of stellar
subdivisions is needed only to show that stellar equivalence classes of stellar
$n$--manifolds are, in fact, the same as piecewise linear homeomorphism classes of combinatorial $n$--manifolds. 

Suppose that
$A$  is any simplex in a simplicial complex 
$K$.  The operation $(A,a)$, of starring $K$ at a point $a$ in the interior of $A$, is the operation 
that changes $K$ to $K^\prime$ by removing $\hbox {\rm st}\, (A,K)$ and replacing it with 
$a\star \partial A\star \hbox {\rm lk}\, (A,K)$.  This is written  $K^\prime = (A,a)K$ or $K \buildrel
{(A,a)} \over
\longmapsto K^\prime$. In the abstract setting $a$ is just a vertex not in $K$.
This type of operation is called a stellar
{\it subdivision}, the inverse operation 
$(A,a)^{-1}$ that changes $K^\prime$ to $K$ is called a  stellar  {\it weld}.  If simplicial complexes $K_1$ and
$K_2$ are related by a sequence of starring operations (subdivisions or welds) and simplicial
isomorphisms, they are called stellar equivalent, written  $K_1 \sim K_2$.
Note that $((A,a)K)\star L = (A,a)(K\star L)$, so that $K_1 \sim K_2$ and $L_1 \sim L_2$ implies that
$K_1\star L_1 \sim K_2\star L_2$.

There now follow the definitions of stellar balls, spheres, and manifolds.  These ideas are not normally
encountered as it will follow, at the end of Section 4, that they are the same as the more familiar 
combinatorial  balls, spheres, and manifolds; in turn these are, up to piecewise linear equivalence, the
even more familiar piecewise linear balls, spheres, and manifolds. 

\proclaim {Definition 3.1}\rm  A stellar $n$--ball is a simplicial complex $B^n$ stellar equivalent to $\Delta ^n$.
A stellar $n$--sphere is a simplicial complex $S^n$ stellar equivalent to $\partial\Delta ^{n+1}$. 
A stellar $n$--manifold is a simplicial complex $M$ such that, for every vertex $v$ of $M$,
$\hbox {\rm lk}\, (v,M)$ is a stellar $(n-1)$--ball or a stellar $(n-1)$--sphere. \endproclaim

A first exercise with stellar ideas is to prove, for stellar balls and spheres, the following:
 $$B^m\star B^n \sim B^{m+n+1};\ S^m\star S^n \sim S^{m+n+1} ;\ B^m\star S^n \sim  B^{m+n+1}.$$

\proclaim {Lemma 3.2}  Suppose $M$ is a stellar $n$--manifold.
\roster
\item If $A$ is an $s$--simplex of $M$ then $\hbox {\rm lk}\, (A,M)$ is a stellar
$(n-s-1)$--ball or $(n-s-1)$--sphere.
\item If $M \sim M^\prime$ then $M^\prime$ is a stellar
$n$--manifold.  
\endroster
\endproclaim

 \demo {Proof}  Assume inductively that both parts of the Lemma are true for stellar $r$--manifolds where $r < n$. 
Part (ii) implies that, for $r < n$, stellar $r$--balls and stellar $r$--spheres are stellar $r$--manifolds. 

Suppose that   $A$ is an $s$--simplex of $M$.  Writing $A = v\star B$, for $v$ a vertex of $A$ and
$B$ the opposite face, $\hbox {\rm lk}\, (A,M) = \hbox {\rm lk}\, (B, \hbox {\rm lk}\, (v,M))$. But $\hbox {\rm lk}\, (v,M)
$ is a stellar
$(n-1)$--sphere or ball and so is a stellar
$(n-1)$--manifold by induction and, again by the induction, $\hbox {\rm lk}\, (B, \hbox {\rm lk}\, (v,M))$ is a 
stellar $(n-s-1)$--ball or
$(n-s-1)$--sphere. This establishes the induction step for (i). 

Suppose that $M_1$ and $M_2$ are complexes and $(A,a)M_1 = M_2$.  If a vertex $v$ of $M_1$ is not in $\hbox {\rm
st}\, (A,M_1)$ then 
$\hbox {\rm lk}\, (v,M_2)  = \hbox {\rm lk}\, (v,M_1)$.  If $v \in \hbox {\rm lk}\, (A,M_1) $, then 
$\hbox {\rm lk}\, (v,M_2) =  (A,a)\hbox {\rm lk}\, (v,M_1)$.  Thus the links of $v$ in $M_1$ and $M_2$ are related by
a stellar move; if one is a stellar ball or sphere, then so is the other.  If $v$ is a vertex of $A$ so that
$A = v\star B$, then
$\hbox {\rm lk}\, (v,M_2)$ is isomorphic to $(B,b)\hbox {\rm lk}\, (v,M_1)$ so similar remarks apply.  Finally   
$\hbox {\rm lk}\, (a,M_2)$ is $\partial A \star \hbox {\rm lk}\, (A,M_1)$.  If $M_1$ is a stellar $n$--manifold,
$\hbox {\rm lk}\, (A,M_1)$ is, by (i), a stellar ball or sphere and hence so is $\partial A \star \hbox {\rm
lk}\, (A,M_1)$. Hence if two manifolds differ by one stellar move and one of them is a stellar $n$
manifold then so is the other. \enddemo

The boundary $\partial M$ of a stellar $n$--manifold $M$ is all simplexes that have as link in $M$ a
stellar ball.  It is not hard to see that this is a subcomplex and is a stellar $(n-1)$--manifold without
boundary.  

If $X$ is a stellar $n$--ball it follows at once from the definitions that $X \sim \Delta^n \sim v \star \partial
\Delta^n
\sim v \star \partial X$. However, in a sequence of stellar moves that produces this equivalence many of the moves may
be of the form
$(A,a)^{\pm 1}$ where 
$A \in \partial X$.  If $A \notin \partial X$, $(A,a)$ is called an {\it internal move}.

\proclaim {Definition 3.3}\rm If $X$ is a stellar $n$--ball and there is an equivalence
$X \sim  v \star \partial X$ using only {\it internal} moves, then $X$ is said to be {\it starrable} and the sequence
of moves is  a starring of $X$.          \endproclaim

\proclaim {Lemma 3.4} Suppose that a stellar $n$--ball $K$ can be starred. Then any stellar
 subdivision $L$ of $K$ can also be starred. \endproclaim

 \demo {Proof}   Suppose that $K \buildrel {(A,a)} \over \longmapsto L$.  It may be
assumed that $A\in \partial K$ for otherwise the result follows trivially. 
Suppose that $K \sim v \star  \partial K $ by $r$ (internal) moves and suppose,
inductively that the result is true for any stellar $n$--ball that can be
starred with less than $r$ moves. If $r=0$ both $K$ and $L$ are cones and there
is nothing to prove.  

Suppose that the first of the $r$ moves is a stellar {\it subdivision} $K \buildrel {(B,b)} \over
\longmapsto K_1$.  As this is an internal move, $\partial K_1 = \partial
K$ and so $A \in \partial K_1$.  Let $L_1$ be the result of starring $K_1$ at
$a$. It follows from the induction hypothesis that $L_1$ can be starred. 

If there
is no simplex in
$K$ with both $A$ and $B$ as a face, or if $A
\cap B = \emptyset$, then $L \buildrel {(B,b)} \over \longmapsto L_1$ so it
follows at once that $L$ can be starred. Note that essentially the proof here is the assertion
$(A,a)(B,b)$ and $(B,b)(A,a)$ both change $K$ to $L_1$.

Now suppose that $A \cap B = C$ for some $C \in K$ and, writing $A = A_0 \star C$ and $B = B_0 \star C$, that
$A_0 \star B_0 \star C \in K$.  Compare again the results of the subdivisions $(A,a)(B,b)$ and $(B,b)(A,a)$ on
$K$.  The resulting complexes can only differ in the way that the join of $A_0 \star B_0 \star C$ to its link is
subdivided.  Thus consider the effect of $(B,b)(A,a)$ on $A_0 \star B_0 \star C$.  
$$ \eqalign{
A_0 &\star B_0 \star C  \buildrel {(A,a)} \over \longmapsto  \partial A \star a \star B_0 \cr
& = (A_0 \star \partial C \star a \star B_0) \;\cup\; (C \star \partial A_0 \star a \star B_0 )\cr
& \buildrel {(B,b)} \over \longmapsto (A_0 \star \partial C \star a \star B_0 ) \;\cup\; (\partial A_0 \star a \star b \star
\partial B) \cr
& = (A_0 \star \partial C \star a \star B_0) \;\cup\; (\partial A_0 \star a \star b \star B_0 \star \partial C) \;\cup\; (\partial
A_0 \star a \star b \star \partial B_0 \star C)
\cr 
& = (\partial (b \star A_0) \star a \star B_0 \star \partial C) \;\cup\; (\partial A_0 \star a \star b \star \partial B_0 \star C)
\ . \cr} $$ 
If $d$ is a point in the interior of $a \star B_0$ the stellar subdivision $(a \star B_0,d)$ changes this last
complex to
$$(\partial (b \star A_0) \star d \star \partial(a \star B_0) \star \partial C) \;\cup\; (\partial A_0 \star a \star b \star
\partial B_0 \star C),$$
see Figure 2.

\figure\relabelbox\small
\epsfxsize0.7\hsize
\epsfbox[37 310 582 441]{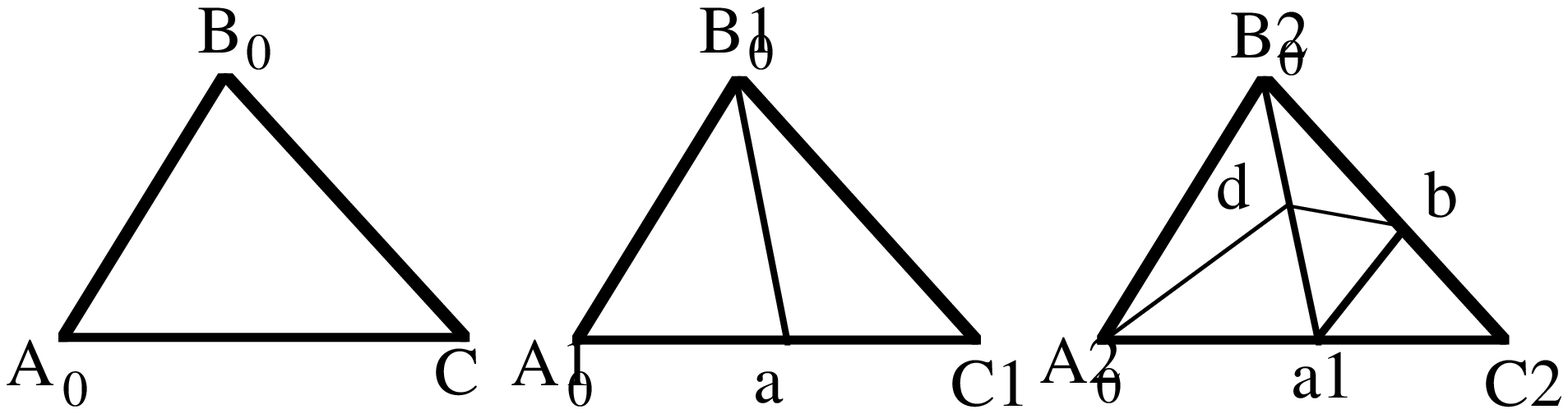}
\relabel {A}{$A_0$}
\relabel {A1}{$A_0$}
\relabel {A2}{$A_0$}
\relabel {B}{$B_0$}
\relabel {B1}{$B_0$}
\relabel {B2}{$B_0$}
\relabel {C}{$C$}
\relabel {C1}{$C$}
\relabel {C2}{$C$}
\relabel {a}{$a$}
\relabel {a1}{$a$}
\relabel {b}{$b$}
\relabel {d}{$d$}
\endrelabelbox
\endfigure

 This is a symmetric expression with respect to $A$ and $B$ so the subdivision 
$(b \star A_0,d^\prime)(A,a)(B,b)$, for $d^\prime$ in the interior of $b \star A_0$, produces  on $A_0 \star B_0 \star C$ an
isomorphic  result.  Taking joins with the link of
$A_0 \star B_0 \star C$ in $K$, it follows that $(b \star A_0,d^\prime)L_1$ and $(a \star B_0,d)(B,b)L$ are isomorphic.  Thus $L$
differs from $L_1$ by internal moves and so is starrable.  

If the first of the $r$ moves is a weld on $K$ that creates $K_1$ then $K_1 \buildrel {(B,b)} \over \longmapsto K$
for some simplex $B$ in the interior of $K_1$. 
However it has just been shown that then $L$ and $L_1$ differ by internal moves so, as $L_1$ is starrable, so is
$L$.  \enddemo
\proclaim {Lemma 3.5} Suppose that a stellar $n$--ball $K$ can be starred. Then $v \star K$, the cone on $K$ with vertex $v$,
can also be starred. \endproclaim

 \demo {Proof}   Suppose that $K \sim u \star (\partial K)$ by $s$ internal moves and suppose,
inductively that the result is true for any stellar $n$--ball that can be
starred with fewer  than $s$ moves.  If $s=0$ then $K = u \star (\partial K)$ and $v \star u \star (\partial K)$ can be starred by a
stellar subdivision at a point in the interior of the 1--simplex $v \star u$. 

Suppose that the first of the $s$ moves changes $K$ to $K_1$ by a weld.  Then $v \star K$ can be obtained from $v \star K_1$ by a
stellar subdivision.  As by induction $v \star K_1$ is starrable, it follows from Lemma 3.4 that $v \star K$ is
starrable.  

Thus suppose that, for $A$ in the interior of $K$, the first of the $s$ moves is 
$K \buildrel {(A,a)} \over \longmapsto K_1$.  Let $P = \hbox {\rm\, lk}(A,K)$ and let $Q$ be the closure of
$K - \hbox {\rm\, st}(A,K)$ so that 
$$v \star K = v \star Q \; \cup \; v \star A \star P.$$  
For $b$ a point in the interior of $v \star A$,
$$v \star A \star P \buildrel {(v \star A,b)} \over \longmapsto b \star \partial (v \star A) \star P = v \star b \star \partial A
\star P
\;\cup\; b \star A \star P .$$
 Now $(v \star Q) \cup (b \star v \star  \partial A \star P)$ is a copy of $v \star K_1$ which is starrable by induction. So for some
vertex $w$
$$v \star K \sim w \star \bigl( v \star \partial K  \;\cup\;  Q  \;\cup\;  b \star \partial A \star P  \bigr) 
  \;\cup\; b \star A \star P $$
by internal moves.  However, replacing $v$ by $w$ in the above argument, 
$$(w \star b \star \partial A \star P) \;\cup\; (b \star A \star P) \sim w \star A \star P$$
 also by internal moves.  Hence
$v \star K \sim w \star \partial (v \star K)$  by internal moves.  \enddemo

\proclaim {Theorem 3.6} Any  stellar $n$--ball $K$ can be starred. \endproclaim

 \demo {Proof}   Assume inductively that stellar $(n - 1)$--balls $K$ can always be starred.

Suppose that $K$ is equivalent to $\Delta ^n$ by $r$  moves and 
suppose  the result is true for any stellar $n$--ball that is equivalent to $\Delta ^n$ by
 fewer  than $r$ moves.  Suppose that the first of the $r$ moves changes $K$ to $K_1$.
The induction step follows immediately if this first move is internal and it follows  from Lemma 3.4
if the first move is a weld.  
Thus suppose that $K \buildrel {(A,a)} \over \longmapsto K_1$ where $A \in \partial K$.  Express
$A$ as $v \star B$, where
$v$ is  a vertex of $A$ and $B$ is the opposite face, so that $\hbox {\rm\, st}(A,K) = v \star B \star \hbox {\rm\, lk}(A,K)$.  
However  $B \star \hbox {\rm\, lk}(A,K)$,
being a stellar $(n-1)$--ball (by Lemma 3.2) is starrable by induction on $n$, and hence $\hbox {\rm\, st}(A,K)$ is
starrable by Lemma 3.5.  Denoting by $X$ the closure of $K - \hbox {\rm\, st}(A,K)$ it follows that, by
internal moves,
$$K \sim X \;\cup\; w  \star \partial \hbox {\rm\, st}(A,K)  
= X \;\cup\;  w \star \bigl( \partial A \star \hbox {\rm\, lk}(A,K) \cup A \star \hbox {\rm\, lk}(A,\partial K)\bigr) .$$
But $X \;\cup\;  w \star \partial A \star \hbox {\rm\, lk}(A,K)$ is isomorphic to $K_1$ and so it is starrable
and can thus be changed by internal moves to 
$x \star \partial\bigl( X \;\cup\;  w \star \partial A \star \hbox {\rm\, lk}(A,K)\bigr) $.  
However, 
$$\bigl( x \star w \star \partial A \star \hbox {\rm\, lk}(A,\partial K)\bigr) 
\;\cup\; \bigl( w \star A \star \hbox {\rm\, lk}(A,\partial K)\bigr) $$
 can be
changed by one internal weld to $x \star A \star \hbox {\rm\, lk}(A,\partial K) = x \star \hbox {\rm\, st}(A,\partial K)$.  Thus $K \sim
x \star \partial K$ by internal moves, as required to complete the induction argument.   \enddemo

\proclaim {Lemma 3.7} Let $M$ be a stellar $n$--manifold containing a stellar $(n-1)$--ball $K$ in its boundary.  Suppose that
the cone $v \star K$ intersects $M$ only in $K$. Then $v \star K \cup M \sim M$. \endproclaim

 \demo {Proof}   By Theorem 3.6 the ball $K$ is starrable.  Any stellar subdivision of $K$ extends to a stellar subdivision 
of
$M$ and of $v \star K$.  An internal weld on $K$ certainly extends over $v \star K$ and it can be extended over $M$, after stellar
moves, in the following way.  Suppose the weld is $K_1 \buildrel {(A,a)} \over \longmapsto K$.  Using Theorem 3.6
again, star the stellar
$n$--ball $\hbox {\rm\, st}(a,M)$, thus changing $M$ to a stellar equivalent manifold $M^\prime$.  This $M^\prime$
contains a cone on $\hbox {\rm\, st}(a,K)$ which can be used to extend the weld over $M^\prime$.  Thus it may be
assumed that $K$ is a cone, $x \star \partial K$  say, and (after starring $\hbox {\rm\, st}(x,M)$) that $M$ contains a cone,
$y \star K$ say, on $K$.  The result then follows from the observation that
$$(v \star x \star \partial K) \,\cup\, (y \star x \star \partial K) \sim v \star y \star \partial K \cong v \star x \star \partial
K .\eqno{\qed}$$

\proclaim {Theorem 3.8} {\rm(Newman \cite 9)}\qua  Suppose $S$ is a stellar $n$--sphere containing a stellar
$n$--ball
$K$.  Then the closure of $S - K$ is a stellar  $n$--ball. \endproclaim

 \demo {Proof}   Inductively assume the theorem is true for all stellar $(n-1)$--spheres.  This assumption implies that
 the closure of $S - K$, is a stellar  $n$--manifold; it is used in checking that the link   of a vertex in the boundary
 is indeed a stellar $(n-1)$--ball.

As, by Theorem 3.6, $K$ can be starred, it is sufficient to prove the result when $K$ is restricted to being the
star of a vertex of $S$.  The sphere $S$ is equivalent to $\partial \Delta^{n+1}$ by a sequence of $r$ stellar
moves.  Assume inductively that the (restricted) result is true for fewer moves.  Suppose that   
$S \buildrel {(A,a)} \over \longmapsto S_1$ is the first of the $r$ moves and that $v$ is a vertex of $S$. 
Let $X$ be the closure of $S - \hbox {\rm\, st}(v,S)$ and $X_1$ be the closure of $S_1 - \hbox {\rm\, st}(v,S_1)$.
If $v \not\in A$ then $X \buildrel {(A,a)} \over \longmapsto X_1$.  If $A = v \star B$ for some face $B$ of $A$, then
$X \cup a \star B \star \hbox {\rm\, lk}(A,S) = X_1$.  Now $X$ is a stellar $n$--manifold containing the stellar $(n-1)$--ball $B \star
\hbox {\rm\, lk}(A,S)$ in its boundary. Thus by Lemma 3.7 $X \sim X_1$.  However induction on $r$ asserts that $X_1$ is a stellar
$n$--ball and so
$X$ is too. 

Suppose the first of the $r$ moves is a weld so that  $S_1 \buildrel {(A,a)} \over \longmapsto S$.  If $v \neq a$
so that $v \in S_1$, then $X \sim X_1$, by the same argument as before, and so $X$ is a stellar $n$--ball.  If $v = a$
then
$X$ is the closure of $S_1 - \hbox {\rm\, st}(A,S_1)$.  Let $x$ be any vertex of $A$, so that $A = x \star B$ say, and let
$Y_1$ be the closure of 
$S_1 - \hbox {\rm\, st}(x,S_1)$. This $Y_1$ is a stellar $n$--ball (induction on $r$) and $\partial Y_1$ contains the
stellar $(n-1)$--ball $B \star \hbox {\rm\, lk}(A,S_1)$.  Let $Z$ be the closure of $\partial Y_1 - B \star \hbox {\rm\, lk}(A,S_1)$, a
stellar $(n-1)$--ball by the induction on $n$.  Then $X = Y_1 \cup x \star Z$ and so $X$ is a stellar $n$--ball by Lemma 3.7.
 \enddemo

\proclaim  {Theorem 3.9}{\rm (Alexander \cite 1)}  Let $M$ be a stellar $n$--manifold, let $J$ be a
stellar $n$--ball. Suppose that 
$M \cap J = \partial M \cap \partial J$ and that this intersection is  a stellar $(n-1)$--ball $K$. 
 Then $M \cup J  \sim M$. \endproclaim

 \demo {Proof}   It may be assumed, by Theorem 3.6, that $J$ is a cone, $x \star \partial J$ say.
Let $L$ be the closure of $\partial J - K$. Then $L$ is a stellar $(n-1)$--ball by Theorem 3.8 and so can
be starred (by Theorem 3.6) to become $y \star \partial L = y \star \partial K$; the stellar equivalence of this starring extends
over the cone
$x \star L$ to give $x \star L \sim x \star y \star \partial L =  y \star x \star \partial K$.  Thus  
$$M \cup J  \sim M \cup x \star K \cup y \star x \star \partial K$$
 and this is equivalent to $M$ by two applications of Lemma 3.7.
 \enddemo


\head  Arbitrary subdivision \endhead

This section connects the idea, just discussed, of stellar moves with that of arbitrary subdivision and, in so doing, 
equates stellar manifolds with combinatorial manifolds. Here it really does help to have simplicial complexes contained in some
${\Bbb R}^N$.
In order to consider arbitrary subdivisions of simplicial complexes it will be convenient to consider the idea  of a convex linear
cell complex, a generalisation of the idea of a simplicial complex.  In principle this is needed  because, if $K_1$ and $K_2$
are simplicial complexes in ${\Bbb R}^N$, there is no natural simplicial triangulation for $ |K_1| \cap |K_2| $.    Any
{$(n-1)$}--dimensional affine subspace (or {\it hyperplane}) ${\Bbb R}^{N-1}$ of ${\Bbb R}^N$ separates ${\Bbb R}^N$ into two
closed half spaces 
${\Bbb R}^N_+$ and ${\Bbb R}^N_-$.  
A convex linear cell (sometimes called a
polytope) in
${\Bbb R}^N$ is a {\it compact} subset that is the intersection of finitely many such half spaces.  A face of the cell is the
intersection of the cell and any subset of the hyperplanes that define these half spaces; the proper faces constitute the
boundary of the cell. The dimension of the cell is the smallest dimension of an affine subspace that contains the cell.  A
convex linear cell complex $C$ in
${\Bbb R}^N$ is  then defined in the same way as is a finite simplicial complex but using convex linear cells in place of simplexes.
In the same way the underlying polyhedron $|C|$ and the idea of a subdivision of one convex linear cell complex by another are
defined. 
 Note that any convex linear cell complex has a subdivision that
is a simplicial complex; an example of a simplicial subdivision is a first derived subdivision, one that subdivides  cells in
some order of increasing dimension, each as a cone on its subdivided boundary.

 \proclaim {Definition 4.1}\rm Two convex linear cell complexes are piecewise linearly homeomorphic if they have
simplicial subdivisions that are isomorphic. \endproclaim

\noindent It is easy to see that if there is a face-preserving bijection between the cells of one convex linear cell complex and
those of another, then the complexes are piecewise linearly homeomorphic.

\proclaim {Lemma 4.2} Let $A$ be a convex linear $n$--cell and let $\alpha\partial A$ be any subdivision of its boundary.
Then $\alpha\partial A$ is piecewise linearly isomorphic to
$\partial \Delta^{n}$. \endproclaim

 \demo {Proof}  Choose $\Delta^{n}$ to be contained (linearly) in $A$ and let $v$ be a point in the interior of $\Delta^{n}$.
Let $\{ B_i\}$ be the cells of $\alpha\partial A$ and $\{ C_j\}$ be the simplexes of $\partial \Delta^{n}$.  Let
$D_{i,j} = (v \star B_i) \cap C_j$.  Then $\{ D_{i,j}\}$ forms a convex linear cell complex  subdividing $\partial \Delta^{n}$.
Let $\beta \partial \Delta^{n}$ be a simplicial subdivision of this cell complex.  Projecting, radially from $v$, the simplexes of 
$\beta \partial \Delta^{n}$ produces a simplicial subdivision $\gamma \alpha\partial A$ of $\alpha\partial A$ and the radial
correspondence gives the required bijection between the simplexes of $\beta \partial \Delta^{n}$ and $\gamma \alpha\partial A$
(though such a projection is not linear when restricted to an actual simplex).  \enddemo

\proclaim {Corollary 4.3} Any subdivision of a convex linear $n$--cell is piecewise linearly isomorphic to 
$\Delta^{n}$. \endproclaim

 \demo {Proof}  Let  $\alpha A$ be a subdivision of a convex linear $n$--cell $A$.  As above, 
$v \star \gamma \alpha\partial A$ is a
simplicial complex isomorphic to the subdivision 
$v \star \beta \partial \Delta^{n}$ of $\Delta^{n}$.  The intersection of the simplexes of  $v \star \gamma \alpha\partial A$ with
the cells of
$\alpha  A$ gives a common cell-subdivision of these two complexes.  \enddemo

The aim of what  follows in Theorem 4.5 is to prove the fact that any simplicial subdivision of a simplicial complex $K$ is
stellar equivalent to $K$. The next lemma is a weak form of this.

\proclaim {Lemma 4.4}  Let $K$ be an $n$ dimensional simplicial complex and let  $\alpha K$ be a simplicial subdivision of
$K$ with the property that, for each simplex $A$ in $K$, the subdivision $\alpha A$ is a stellar ball. 
Then $\alpha K \sim K$. \endproclaim

 \demo {Proof}   Let $\beta_rK$ be the subdivision of $K$ such that, if $A \in K$ and  $\hbox {\rm dim}\,
A \leq r$, then 
$\beta_rA = \alpha A$ and if $\hbox {\rm dim}\, A > r$ then $A$ is subdivided as  the cone on the already defined
subdivision of its boundary (think of the simplexes being subdivided one by one in some order of increasing
dimension).  If $A \in K$ and $\hbox {\rm dim}\, A = r$ then the stellar $r$--ball $\alpha A$ can be starred.  As this is
an equivalence by internal moves, simplexes of dimension less than $r$ are unchanged, and the stellar moves can be
extended conewise over the subdivided simplexes of higher dimension.  In this way,  $\beta_rK \sim \beta_{r-1}K$.
However $\beta_nK = \alpha K$ and $\beta_0K$ is just a first derived subdivision of $K$ which is certainly
stellar equivalent to $K$. \enddemo

The next theorem is the promised result that links arbitrary subdivisions with stellar moves.  It is an easy exercise to
produce a subdivision of $\Delta^2$ that is not the result of a sequence of stellar {\it subdivisions} on  $\Delta^2$. 
It is a classic conjecture, which the author believes to be still unsolved, that if two complexes are piecewise linearly
homeomorphic (that is, they have isomorphic subdivisions) then they have isomorphic subdivisions each obtained by a sequence of
stellar {\it subdivisions} (with no welds being used).

\proclaim {Theorem 4.5} Two $n$--dimensional simplicial complexes are piecewise linearly homeomorphic if and only if they are
stellar equivalent. \endproclaim

 \demo {Proof}   Clearly, stellar equivalent complexes are piecewise linearly homeomorphic.  Thus it is sufficient to prove
that if $\alpha K$ is a simplicial subdivision of any $n$--dimensional simplicial complex 
$K$, then $\alpha K \sim K$. Assume inductively that this is true for all simplicial complexes of dimension less
than $n$.

Suppose that $C$ is a convex linear cell complex contained in ${\Bbb R}^N$,  that ${\Bbb R}^{N-1}$ is an
affine $(N-1)$--dimensional subspace of ${\Bbb R}^N$ and that ${\Bbb R}^N_+$ and ${\Bbb R}^N_-$ are the closures
of the complementary domains of ${\Bbb R}^{N-1}$.  The {\sl slice of $C$ by ${\Bbb R}^{N-1}$} is defined to be 
the convex linear cell complex consisting of all cells of one of the forms  $A \cap {\Bbb R}^N_+$, 
$A \cap {\Bbb R}^N_-$ or $A \cap {\Bbb R}^{N-1}$ where $A$ is a cell of $C$.  An $r$--slice subdivision of $C$ is the
result of a sequence of $r$ slicings by  such  affine subspaces of dimension $(N-1)$.
Let $K$ be an $n$--dimensional simplicial complex and suppose $C_r$ is an $r$--slice subdivision of $K$.  Let $\beta
C_r$ be any simplicial subdivision of $C_r$ such that each top dimensional ($n$--dimensional) cell is subdivided as
a cone on the subdivision of its boundary.  Suppose inductively that for smaller $r$ any such $\beta C_r$ is
stellar equivalent to $K$.  

To start the induction, when $r=0$, it is necessary to consider a simplicial subdivision $\beta K$ of $K$ in which all the
$n$--simplexes of $K$ are subdivided as cones on their boundaries.  If $A \in K$ and $\hbox {\rm dim}\, A < n$ then 
 $\beta A$ is, by the induction on $n$, a stellar ball.  If $\hbox {\rm dim}\, A = n$, then $\beta A$ is a cone on 
$\beta\partial A$.  But $\beta\partial A$ is, by the induction on $n$, a stellar $(n-1)$--sphere
and so the cone $\beta A$ is a stellar $n$--ball.   Thus in $\beta K$ the subdivision of every simplex of $K$ is a stellar
ball and so, by Lemma 4.4, $K \sim \beta K$. 

Suppose that the $r$--slice subdivision $C_r$ of $K$ is obtained by slicing the $(r-1)$--slice subdivision $C_{r-1}$.  Choose
a simplicial subdivision $\gamma$ of $C_{r-1}$ so that, on the cells of dimension less than $n$, $\gamma$ is the the
slicing by the $r^{\hbox {\rm th}}$ affine subspace ${\Bbb R}^{N-1}$ followed by $\beta$. On an $n$--cell,  $\gamma$ is a
subdivision of the cell as the cone on its subdivided boundary. The induction on $r$ implies that $\gamma C_{r-1} \sim K$. 
Thus, to complete the induction it is necessary to show that $\gamma C_{r-1} \sim \beta C_r$. Now 
$\gamma C_{r-1}$ and $\beta C_r$ differ only in the way that the interiors of the $n$--cells of $C_{r-1}$ are
subdivided.  Suppose that $A$ is an $n$--cell of $C_{r-1}$ and that $A \cap {\Bbb R}^{N-1}$ is an $(n-1)$--cell.
Thus in $C_r$, the cell $A$ is divided into two cells, $X$ and $Y$ say, intersecting in the cell $A \cap {\Bbb R}^{N-1}$
contained in their boundaries.  Now $\beta\partial X$ is, by induction on $n$ and Lemma 4.2, a stellar $(n-1)$--sphere
and so
$\beta X$ , which is  the cone on this, is a stellar
$n$--ball. Similarly $\beta Y$ is a stellar
$n$--ball.
Again by induction on $n$ and using Corollary 4.3, $\beta(A \cap {\Bbb R}^{N-1})$ is a
stellar $(n-1)$--ball.  Thus by Theorem 3.9 $\beta X \cup \beta Y$ is a stellar $n$--ball and so, by Theorem 3.6 can be starred by
internal moves.  This starring process changes $\beta A$ to $\gamma A$.  Repetition of this on every $n$--cell  of
$C_{r-1}$ shows that $\gamma C_{r-1} \sim \beta C_r$.  That completes the proof of the induction on $r$.

\figure\relabelbox\small
\epsfxsize0.57\hsize
\epsfbox{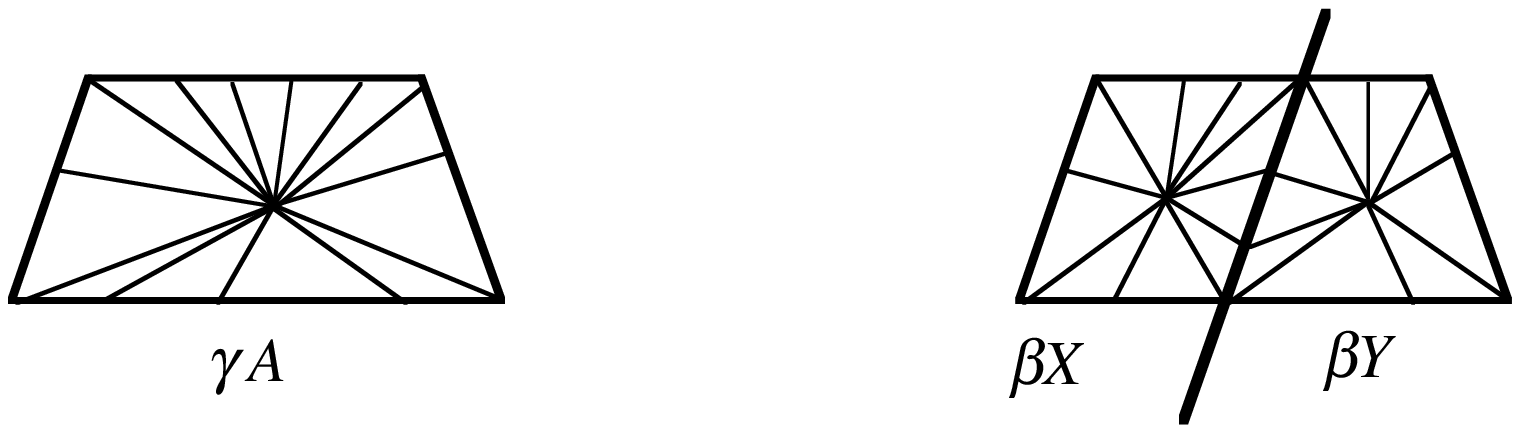}
\adjustrelabel <-5pt,0pt> {A}{$\gamma A$}
\adjustrelabel <-2pt,0pt> {X}{$\beta X$}
\relabel {Y}{$\beta Y$}
\endrelabelbox\endfigure

Finally, consider the simplicial  subdivision $\alpha K$ of the   $n$ dimensional simplicial complex $K$ assumed to
be contained in some  ${\Bbb R}^N$.  For each simplex $A$ of $\alpha K$ (of dimension less than $N$)
choose an affine subspace ${\Bbb R}^{N-1}$ containing $A$ but otherwise in general position with respect to all the vertices of
$\alpha K$.  Using all these affine subspaces,   construct an $r$--slice subdivision $C$ of $K$. 
Suppose $A\in \alpha K$ is a $p$--simplex contained in a $p$--simplex $B$ of $K$. The above general position requirement ensures
that
$B$ meets the copy of ${\Bbb R}^{N-1}$, selected to contain  any $(p-1)$--dimensional face of $A$, in the intersection of $B$ with
the $(p-1)$--dimensional affine subspace containing that face (and not in the whole of $B$). Thus the slice subdivision of $K$
using just these particular $p+1$ hyperplanes contains the simplex $A$ as a cell.  Of course, $A$ is subdivided if more of
the slicing operations are used, but this means that
$C$ is also an  $r$--slice subdivision of $\alpha K$.  Let $\beta$ be a simplicial subdivision of $C$ as described above. 
Then by the above $K \sim \beta C \sim \alpha K$.  \enddemo

In the light of the preceding theorem, combinatorial balls, spheres and manifolds are the same as
the stellar balls, spheres and manifolds previously considered. Thus in all the preceding results of section 3, the word `combinatorial' can now
replace `stellar'.  In the rest of this paper the `combinatorial' terminology will be used.
It is worth remarking that in \cite {10} Newman modifies Theorem 3.6 so that he can improve Theorem 4.5 by restricting the
needed elementary stellar subdivisions and welds to those involving only 1--simplexes. 

As in Lemma 3.2, it follows that the link of any
$r$--simplex in
$M$ is a combinatorial $(n-r-1)$--ball or sphere. 
The following is a short technical lemma concerning
this that will be used later.

\proclaim {Lemma 4.6} Suppose that $B$ is an $r$ simplex, that $X$ is any simplicial complex and $ \partial B \star X$
is a combinatorial $n$--sphere or $n$--ball.  Then $X$ is a combinatorial $(n-r)$--sphere or $(n-r)$--ball. \endproclaim

 \demo {Proof}   For $v$ a new vertex, $v \star \partial B \star X$ is a combinatorial $(n+1)$--ball.  But this is
piecewise linearly homeomorphic to $B \star X$. Thus
$\hbox {\rm lk}\, (B,B \star X)$, namely
$X$, is a  combinatorial $(n-r)$--sphere or $(n-r)$--ball.  \enddemo


\head  Moves on manifolds 
\endhead

This final section deduces that bistellar moves suffice to change one triangulation of a closed
piecewise linear manifold to another.  A similar theorem for bounded manifolds is also included. 
A version of the proof of Pachner \cite {11}, \cite {12}, is employed. Firstly there   follow 
definitions of  some relations between  combinatorial
$n$--manifolds. They are best thought of as more types of `moves' changing one complex to another within the same  piecewise
linear homeomorphism class.  

\proclaim {Definition 5.1}\rm Let $A$ be a non-empty simplex in a combinatorial $n$--manifold $M$ such that
$\hbox {\rm lk}\, (A,M) = \partial B \star L$ for some simplex $B$ with $\emptyset \neq B \notin M$ and some
complex $L$. Then $M$ is related to $M^\prime$ by the stellar exchange $\kappa (A,B)$, written 
$M \buildrel \kappa (A,B) \over  \longmapsto   M^\prime$, if  $M^\prime$ is obtained by removing 
$A \star \partial B \star L$ from $M$ and inserting $\partial A \star B \star L$.   \endproclaim

\figure\relabelbox\small
\epsfxsize0.52\hsize
\epsfbox[124 310 500 459]{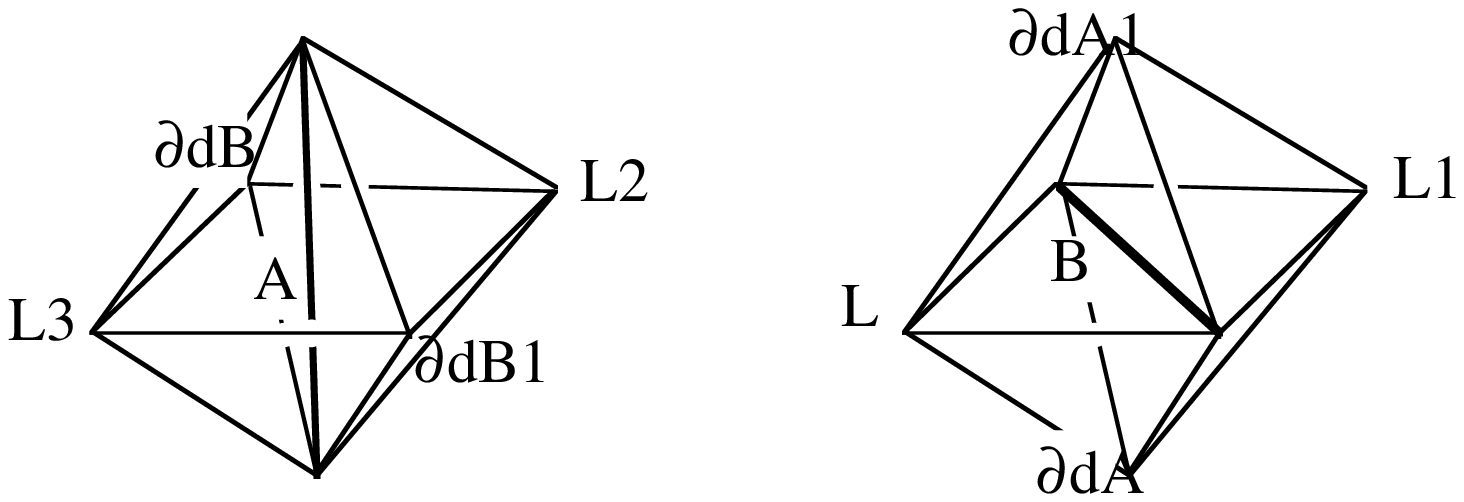}
\adjustrelabel <-2pt, 0pt> {A}{$A$}
\adjustrelabel <-1pt,-1pt> {B}{$B$}
\adjustrelabel <-7pt, 0pt> {dA}{$\partial A$}
\adjustrelabel <-6pt, 0pt> {dB}{$\partial B$}
\adjustrelabel <0pt, 5pt> {dA1}{$\partial A$}
\relabel {dB1}{$\partial B$}
\adjustrelabel <-2pt, 0pt> {L}{$L$}
\adjustrelabel <-2pt, 0pt> {L1}{$L$}
\relabel {L2}{$L$}
\relabel {L3}{$L$}
\endrelabelbox
\endfigure

This idea is illustrated in Figure 4.  If $L = \emptyset$, then $\kappa (A,B)$ is, as defined in
Section 2,  a bistellar move. As in Section 2, for this definition and for others like it, it is
best to consider
$M$ as an abstract simplicial complex. Otherwise, topologically,    `$B \notin M$' should here be
expanded to mean   that $A \star B \star L$ exists and meets $M$
only in $A \star \partial B \star L$. It is important here not to be deterred by a feeling that the existence of
simplexes such as $A$ is rather unlikely. Note that $M^\prime \buildrel \kappa (B,A) \over
 \longmapsto   M$ is the inverse move to $\kappa (A,B)$.  If $B$ is a single vertex
$a$, so that $\partial B = \emptyset$, then  $\kappa (A,B)$ is the stellar subdivision $(A,a)$
discussed at length in Section 3. 
Note  that any
$\kappa (A,B)$ 
is the composition of a stellar subdivision and a weld, namely 
$(B,a)^{-1}(A,a)$.     

\proclaim {Definition 5.2}\rm Suppose that $A$ and $B$ are simplexes of a combinatorial $n$--manifold  $M$ with boundary
$\partial M$, that the join
$A \star B$ is an
$n$--simplex of
$M$,  that 
$A \cap \partial M = \partial A$ and $B \star \partial A \subset \partial M$.  The manifold
$M^\prime$ obtained from $M$  by elementary shelling from $B$ is the closure of $M - (A \star B)$.  \endproclaim

Here taking the closure means adding on the smallest number of simplexes (in this case those of
$A \star \partial B$) to achieve a simplicial complex. The relation between $M$ and $M^\prime$
will be denoted
$M \buildrel (\hbox {\rm sh}\, B)  \over  \longmapsto  M^\prime$. Note that 
$\partial M^\prime$ and $\partial M$ are related by a bistellar move.

\proclaim {Definition 5.3}\rm A combinatorial $n$--ball is shellable if it can be reduced to a single $n$--simplex
by a sequence of elementary shellings. A combinatorial $n$--sphere is shellable if removing some
$n$--simplex from it produces a shellable combinatorial
$n$--ball. \endproclaim

\noindent Note that there are well known examples of combinatorial $n$--spheres and combinatorial $n$--balls that
are {\it not} shellable (see, for example, \cite {14}, \cite 3 or \cite 7).  For future use, three
straightforward lemmas concerning shellability now follow.

\proclaim {Lemma 5.4}  If $X$ is a shellable combinatorial $n$--ball or $n$--sphere then the cone $v \star X$ is
shellable. \endproclaim

 \demo {Proof}  If  $X$ is a  combinatorial $n$--ball, for every elementary shelling $(\hbox {\rm sh}\, B)$  in a shelling sequence for
$X$,   perform the shelling
$(\hbox {\rm sh}\, B)$ on $v \star X$ (with $v \star A$ in place of $A$). This shows that $v \star X$ is shellable. Suppose then that  
$X$ is a  combinatorial
$n$--sphere  and $X - C$ is a shellable ball for some $n$--simplex
$C$.  Then $v \star X \buildrel (\hbox {\rm sh}\, C) \over  \longmapsto  v \star (X - C)$ and the result follows from
the preceding case.  \enddemo

\proclaim {Lemma 5.5}  If $X$ is a shellable combinatorial $n$--ball or $n$--sphere then  
$\partial\Delta^r \star X$ is shellable. \endproclaim

 \demo {Proof}  Suppose  $X$ is a  combinatorial $n$--ball. Work by induction on $r$; when $r = 0$ there is
nothing to prove. Write
$\Delta^r = v \star C$ for some vertex $v$ and $(r-1)$--simplex $C$.  Suppose that the shellings 
$(\hbox {\rm sh}\, B_1), (\hbox {\rm sh}\, B_2), \dots ,(\hbox {\rm sh}\, B_s)$ on $X$ change $X$ to a single
$n$--simplex. Then shellings 
$(\hbox {\rm sh}\, (C \star B_1)), (\hbox {\rm sh}\, (C \star B_2)), \dots ,(\hbox {\rm sh}\, (C \star B_s))$ followed by $(\hbox {\rm sh}\, C)$ change 
$\partial\Delta^r \star X$ to
$v \star \partial C \star X$ and the result follows by induction on $r$ and the last lemma.  Now let $X$ be a
combinatorial
$n$--sphere with $X - D$ shellable for some
$n$--simplex
$D$.  Suppose that $\partial C$ can be reduced to a single $(r-2)$--simplex by the removal of an $(r-2)$--simplex
$C_0$ followed by shellings 
$(\hbox {\rm sh}\, C_1), (\hbox {\rm sh}\, C_2), \dots ,(\hbox {\rm sh}\, C_t)$.  Then shellings
$(\hbox {\rm sh}\, (D \star C_0)), (\hbox {\rm sh}\, (D \star C_1)), \dots ,(\hbox {\rm sh}\, (D \star C_t))$ followed by 
$\hbox {\rm sh}\, (D)$ change $\partial (v \star C) \star X - C \star D$ to 
\hbox{$\partial (v \star C) \star (X - D)$}
which
is shellable by the previous case.   \enddemo

\proclaim {Definition 5.6}\rm  Combinatorial $n$--manifolds $M_1$ and $M_2$ are bistellar equivalent, written 
$M_1 \approx M_2$, if they are related by a sequence of bistellar moves and simplicial isomorphisms. \endproclaim

\proclaim {Lemma 5.7}  If $X$ is a shellable combinatorial $n$--ball, then the cone $v \star \partial X \approx X$.
\endproclaim  

 \demo {Proof}  The proof is by induction on  the number, $r$,  of $n$--simplexes in $X$.  If $r=1$ a single starring
of the simplex $X$ produces
$v \star \partial X$.  Suppose that the first elementary shelling of $X$ is 
$X \buildrel (\hbox {\rm sh}\, B) \over  \longmapsto  X_1$, where 
$A \star B$ is an $n$--simplex of $X$, 
$A \cap \partial X = \partial A$ and $B \star \partial A \subset \partial X$. By the induction on
$r$, $X_1 \approx v \star \partial X_1$.  However, $v \star \partial X_1 \cup A \star B$ is
changed to $v \star \partial X$ by the {  bistellar} move $\kappa (A, v \star B)$.  \enddemo 

\proclaim {Corollary 5.8}  Let $M$ be a  combinatorial $n$--manifold, let $A \in (M - \partial M)$ and suppose that
$\hbox {\rm lk}\, (A,M)$ is shellable.  Then $M \approx (A,a)M$ for $a$ a vertex not in $M$. \endproclaim

 \demo {Proof}   Iteration of Lemma 5.4 shows that $A \star \hbox {\rm lk}\, (A,M)$ is
shellable and so Lemma 5.7 implies that $A \star \hbox {\rm lk}\, (A,M)$ is bistellar
equivalent to the cone on its boundary.  However this cone is
$a \star \partial A \star \hbox {\rm lk}\, (A,M)$.  \enddemo

\proclaim {Theorem 5.9}{\rm  (Newman \cite 9, Pachner \cite {11}, 
\cite {12})}\qua Two closed
combinatorial $n$--manifolds are piecewise linearly homeomorphic if and only if they are bistellar
equivalent. \endproclaim

 \demo {Proof}   By Theorem 4.5 it is sufficient to prove that, if closed combinatorial $n$--manifolds
$M$ and
$M^\prime$ are related by a  stellar
exchange, then they are bistellar equivalent.  
Thus suppose that $\hbox {\rm lk}\, (A,M) = \partial B \star L$  and suppose 
$M \buildrel \kappa (A,B) \over \longmapsto   M^\prime$.  
Suppose that  $L$ can be expressed as 
$L = L^\prime  \star {\Cal S}$ 
where ${\Cal S}$ is some join of copies of $\partial \Delta^p$ for any assorted values of $p$ 
(though ${\Cal S}$ could well be empty).  By Lemma 4.6, $L$ and $L^\prime$
are (possibly empty) combinatorial spheres.  Let $\hbox {\rm dim}\, L^\prime = m$.  
As $L^\prime$ is a stellar $m$--sphere it is related to
$\partial \Delta^{m+1}$ by a sequence of $r$, say, stellar moves.
The proof proceeds by induction on the pair $(m, r)$, assuming the result is true for smaller $m$, or the same
$m$ but smaller $r$, smaller that is than the values under consideration. 
Note that if $m = 0$ the
$0$--sphere
$L^\prime$ can be absorbed into ${\Cal S}$.   For any value of $m$, if $r=0$ then $L^\prime$ is
a copy of $\partial \Delta^{m+1}$ and can again be absorbed into ${\Cal S}$.

However, if $L^\prime$ is empty ($m = -1$), then $L$ is shellable (by Lemma 5.5) and
so (again using Lemma 5.5) $\hbox {\rm lk}\, (A,M)$ and $\hbox {\rm lk}\, (B,M^\prime)$ are both shellable.  Thus, by
Corollary 5.8, 
$$M \approx (A,a)M  = (B,a)M^\prime \approx M^\prime.$$
These remarks start the induction.  

Suppose that
the first of the
$r$ moves that relate  $L^\prime$ to $\partial \Delta^{m+1}$ is the stellar exchange $\kappa (C,D)$
(it is convenient to use the stellar exchange idea to avoid distinguishing subdivisions and
welds).  Thus $C \in L^\prime$ and $\hbox {\rm lk}\, (C,L^\prime ) = \partial D \star L^{\prime\prime}$ for
some complex $L^{\prime\prime}$.

There are two cases to consider, the first being when $D \notin M$. Consider the combinatorial
manifolds $M_1$ and $M_2$ obtained from $M$ in the following way
$$M \buildrel \kappa (A \star C,D) \over  \longmapsto   M_1 \buildrel \kappa (A,B) \over
 \longmapsto   M_2   .$$
Because $\hbox {\rm lk}\, (A \star C,M) = \partial B \star {\Cal S} \star \partial D \star L^{\prime\prime}$ and
$\hbox {\rm dim}\, L^{\prime\prime} < m$, it follows by induction on $m$ (regarding 
$\partial B \star {\Cal S}$ as the `new ${\Cal S}$') that $M \approx M_1$. 
Because $\hbox {\rm lk}\, (A,M_1) = \partial B \star {\Cal S} \star \kappa (C,D)L^{\prime}$
 it follows that $M_1 \approx M_2$   by induction on $r$.

The same $M_2$ can also be achieved in an alternative manner:
$$M \buildrel \kappa (A,B) \over  \longmapsto   M^\prime \buildrel \kappa (B \star C,D) \over
 \longmapsto   M_2 .$$
That this is indeed the same $M_2$ can be inferred from considerations of the symmetry between $A$ and $B$,
but it is also easy to check that 
$A \star C \star \partial B \star \partial D \star {\Cal S} \star L^{\prime\prime}$  
is changed by each pair of moves 
to
$$\Bigl( (B \star D \star \partial A \star \partial C) \cup (C \star D \star \partial A \star \partial B) \Bigr)  \star 
{\Cal S} \star L^{\prime\prime}$$  
and clearly the
changes produced on the rest of $M$ are the same.  However, \break
{$\hbox {\rm lk}\, (B \star C,M^\prime ) = \partial A \star {\Cal S} \star \partial D \star
L^{\prime\prime}$}   so, here again by induction
on $m$,   $M^\prime \approx M_2$.
 Hence $M \approx M^\prime$.  

The second case is when $D \in M$.  A trick reduces this to the first case in the
following way.   It may be assumed that $\hbox {\rm dim}\, D \geq 1$ because if $D$ is just a vertex an
alternative new vertex can be used instead. Write $D = u \star E$ for $u$ a vertex and $E$ the opposite
face in $D$.  Let $v$ be a new vertex not in $M$. 
Consider the manifold $\hat M^\prime$ obtained by
$$M \buildrel \kappa (A \star u,v) \over  \longmapsto   \hat M \buildrel \kappa (A,B) \over
 \longmapsto   \hat M^\prime  $$
or obtained alternatively, as seen by symmetry considerations, by 
$$M \buildrel \kappa (A,B) \over  \longmapsto   M^\prime \buildrel \kappa (B \star u,v) \over
 \longmapsto   \hat M^\prime .$$ 
Because $\hbox {\rm lk}\, (A,\hat M) = \partial B \star \kappa(u,v)L^\prime \star {\Cal S}$ and 
$\kappa(u,v)L^\prime$ is just a copy of $L^\prime$ with $v$ replacing $u$, it follows by the
first case that $\hat M \approx \hat M^\prime$.  
However, because  $\hbox {\rm lk}\, (A \star u,M) =  \partial B \star \hbox {\rm lk}\, (u,L^\prime) \star {\Cal S}$,
the induction on $m$ gives $M \approx \hat M$. Similarly $M^\prime \approx \hat M^\prime$.  
Hence it follows again that $M \approx  M^\prime$.
 \enddemo


An elegant version for bounded manifolds,  of this last result, will now be discussed.  The
proof again is based on Pachner's work. 
 Only an outline of the proof, which runs along the same lines as that of Theorem 5.9, will be
given here.  An  obvious extension of terminology will be used.   If the manifold
$M^\prime$ is obtained from $M$  by an elementary shelling then $M$ will be said to be obtained
from $M^\prime$ by an elementary inverse shelling.  

\proclaim {Theorem 5.10}{\rm  (Newman \cite 9, Pachner \cite {12})}\qua  Two connected combinatorial
$n$--manifolds with non-empty boundary are piecewise linearly homeomorphic if and only if they are
related by a sequence of elementary shellings, inverse shellings and a simplicial isomorphism.
\endproclaim

\demo {Outline proof}  
Suppose $M$ is a connected combinatorial $n$--manifold with $\partial M \neq \emptyset$.  Firstly,
note in the following way that a bistellar move on $M$ can be expressed as a sequence of elementary
shellings and inverse shellings.  In  such a bistellar move $\kappa (A,B)$, the  simplex $A$ is in
the interior of $M$ because $\hbox {\rm lk}\, (A,M) = \partial B$. Find a chain of $n$--simplexes, each
meeting the next in an $(n-1)$--face, that connects $\hbox {\rm st}\, (A,M)$ to an
$(n-1)$--simplex $F$ in
$\partial M$ (see Figure 5(i)).  By a careful use of inverse shellings, add  to the closure of
$\partial M - F$ a collar, a copy of $(\overline{\partial M - F}) \times I$,  see Figure 5(ii). 
Then, by elementary shellings of the union of $M$ and this collar, remove all the $n$--simplexes
of the chain together with $\hbox {\rm st}\, (A,M)$.  Next, using inverse elementary shellings, replace 
$\hbox {\rm st}\, (A,M)$ with $B \star \partial A$ and reinsert the chain of  $n$--simplexes.  Finally
remove the collar by elementary shellings.  The use of the collar   ensures that the
chain of $n$--simplexes, which might meet $\partial M$ in an awkward manner, can indeed be
removed by shellings.

\figure\relabelbox\small\def\ss{\scriptstyle}
\epsfxsize0.65\hsize
\epsfbox[140 345 475 455]{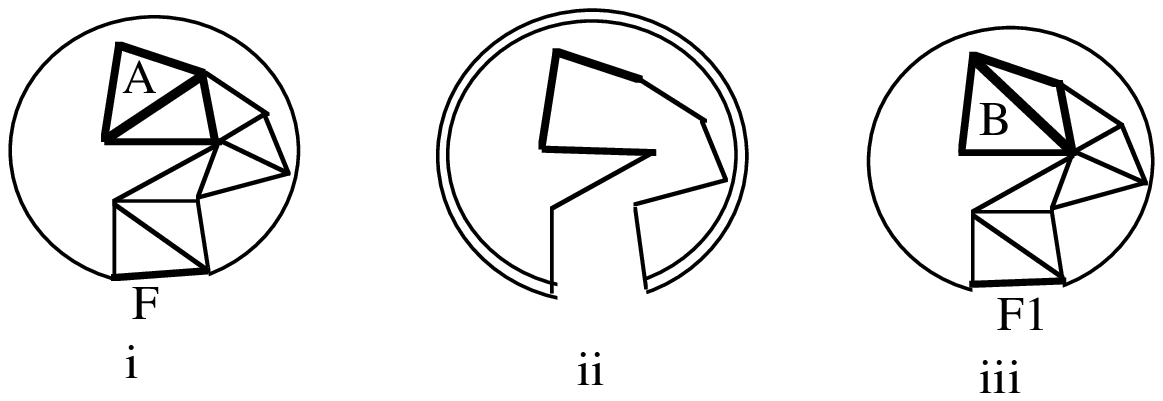}
\relabel {F}{$\ss F$}
\relabel {F1}{$\ss F$}
\adjustrelabel <-1.5pt,-1pt> {A}{$\ss A$}
\adjustrelabel <-1pt,-1pt> {B}{$\ss B$}
\relabel {i}{$\ss\rm(i)$}
\relabel {ii}{$\ss\rm(ii)$}
\relabel {iii}{$\ss\rm(iii)$}
\endrelabelbox
\endfigure

In the light of this last remark and Theorem 4.5, it is sufficient to prove that, if 
\hbox {$M \buildrel \kappa (A,B) \over \longmapsto   M^\prime$} is a stellar exchange, then $M$
and
$M^\prime$ are related by a sequence of elementary shellings and inverse shellings and bistellar
moves.  Proceed as in the proof of Theorem 5.9.  If $A$ is in the interior of $M$, then 
$M$ and $M^\prime$ are related by a sequence of  bistellar moves exactly as in the proof of
Theorem 5.9. If however $A \subset \partial M$ then the proof must be adapted in the following
way. 

Suppose that  $\kappa (A,B)$ is a stellar exchange, where $A \subset \partial M$ and 
$\hbox {\rm lk}\, (A,M) = \partial B \star L$.  As  $\partial B \star L$  is a combinatorial ball
 it follows from Lemma 4.6 that $L$ is a combinatorial ball and $\partial B \subset \partial M$. 
Thus $A \star\partial B \star \partial L \subset \partial M$.   Express $L$ as $L = L^\prime \star
{\Cal S}$ where
${\Cal S}$ is some join of copies of 
$\partial \Delta ^p$ for any assorted values of $p$ and possibly some $\Delta ^q$. 
Lemma 4.6 shows that this $L^\prime$ is a combinatorial sphere or ball.  Note that
${\Cal S}$ and
$\partial {\Cal S}$ are shellable.  Let
$\hbox {\rm dim}\, L^\prime = m$ and suppose that $L^\prime $ is equivalent to $\Delta ^m$ or
$\partial \Delta ^{m+1}$ by way of $r$ stellar moves. If $m = 0$ or $r = 0$, then  $L^\prime $
can be incorporated into ${\Cal S}$ thus reducing $L^\prime $ to the empty set.

If $L^\prime = \emptyset $, whether ${\Cal S}$ is a sphere or a ball,
the ball $L$ and its boundary are both shellable.  In this circumstance the theorem will be 
established below, and that is the start of the inductive proof on the pair $(m, r)$.  
Continue then as in the proof of Theorem 5.9.  in which $\kappa (C,D)$ is a stellar exchange on
$L^\prime$, being the first of the above mentioned $r$ moves.  The proof proceeds exactly as
before except that, if $C$ is in the interior of $L^\prime$, then $A \star C$ and $B \star C$ are in
the interior of $M$ and so the stellar exchanges $\kappa (A \star C, D)$ and $\kappa (B \star C, D)$
are already known to be expressible as a sequence of bistellar moves.  

Finally, for the start of the induction,  consider the situation when $A \subset \partial M$ and
$L$ and $\partial L$, and hence both of
$\hbox {\rm lk}\, (A,M)$ and $\partial \hbox {\rm lk}\, (A,M) = \hbox {\rm lk}\, (A,\partial M)$, are 
shellable.  This means (by Lemma 5.4) that $\hbox {\rm st}\, (A,\partial M)$ is shellable and, using
this, it is straightforward to add to $M$ a cone on  $\hbox {\rm st}\, (A,\partial M)$ by means of a
sequence of elementary inverse shellings.  Let the resulting manifold be $M_+$ (see Figure 6) so that
$M_+ = (v \star \hbox {\rm st}\, (A,\partial M) ) \cup M$ where $v$ is a new vertex.

\figure\relabelbox\small\let\ss\scriptstyle
\epsfxsize0.4\hsize
\epsfbox[219 375 420 488]{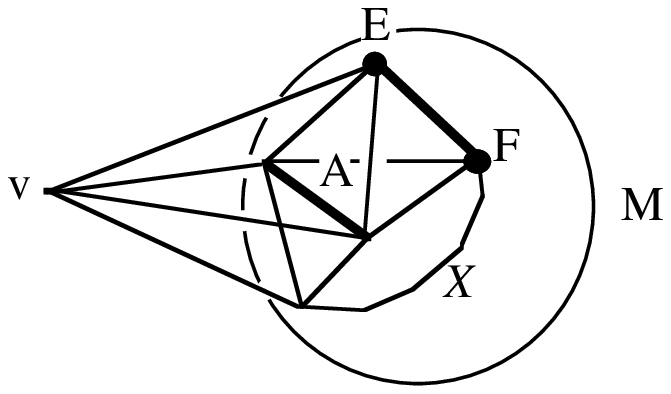}
\adjustrelabel <-2pt, 0pt> {A}{$\ss A$}
\adjustrelabel <-1pt, 0pt> {v}{$\ss v$}
\adjustrelabel <-1pt, 0pt> {E}{$\ss E$}
\relabel {F}{$\ss F$}
\relabel {X}{$\ss X$}
\adjustrelabel <-2pt, 0pt> {M}{$\ss M$}
\endrelabelbox
\endfigure

\noindent For brevity let 
$X = \hbox {\rm lk}\, (A,M)$.
Suppose that the first elementary shelling, in a shelling sequence for $X$, changes $X$ to
$X_1$ by the removal of a simplex $E \star F$, where $F \cap \partial X = \partial F$
$(E \star F) \cap \partial X = E \star \partial F$. Observe that 
$$\hbox {\rm lk}\, (A \star E,M_+) = F \cup (v \star \partial F) = \partial (v \star F) .$$
Thus a stellar exchange $\kappa (A \star E, v \star F)$ can be performed on $M_+$ which has the
effect of removing $\hbox {\rm st}(A \star E, M_+)$ and replacing it with 
$$v \star F \star \partial (A \star E) = v \star F \star ((\partial A \star E) \cup (A \star
\partial E)) = (v \star A \star \partial E \star F) \cup (v \star \partial A \star  E \star F) .$$
The second term is the join of $v \star \partial A$ to the simplex removed from $X$ to create
$X_1$; the first term is the cone from $v$ on $A$ joined to $\partial X_1 - \partial X$.  This
process can now be repeated by using the remaining elementary shellings that change $X_1$ to a single
simplex and then by using the removal of that final simplex as one last `elementary shelling'. The
result is that $M_+$ is bistellar equivalent to 
$$(M - \hbox {\rm st}\, (A, M)) \cup (v \star \partial A \star X)  = (A, v)M.$$
As previously remarked, if $M$ and $M^\prime$ are related by the stellar exchange {$M \buildrel \kappa
(A,B)
\over
\longmapsto   M^\prime$}, then $(A, v)M = (B, v)M^\prime$. 
Of course, $\hbox {\rm lk}\, (B, M^\prime) = \partial A \star L$.
Thus as $(A, v)$ and, similarly, $(B, v)$ have
been shown to be expressible as a sequence of elementary shellings, inverse shellings and bistellar
moves, the combinatorial manifolds $M$ and $M^\prime$ are  related by such moves.  \enddemo

The results of Theorems 5.9 and 5.10 have a simple memorable elegance, albeit the proofs here given have, at times, been a little
involved. It is hoped that publicising these results, together with their proofs, will enable them to find new applications
based on confidence in their veracity.

\Addresses\recd

\enddocument